\documentclass[12pt]{article}

\usepackage[latin1]{inputenc}
\usepackage{amsmath,amssymb}
\usepackage{amsmath,accents}
\usepackage{alltt,graphics,graphpap,color}
\usepackage[dvips]{graphicx}
\usepackage{amsfonts}
\usepackage{amscd}
\usepackage{mathrsfs}

\newtheorem{Theorem}{Theorem}[section]

\newtheorem{Proposition}[Theorem]{Proposition}
\newtheorem{Lemma}[Theorem]{Lemma}
\newtheorem{Corollary}[Theorem]{Corollary}

\newtheorem{Notation}[Theorem]{Notation}

\newcommand{\ep}{\varepsilon}
\newcommand{\ld}{\lambda}
\newcommand{\mm}{{\cal M}}
\newcommand{\bmm}{{\overline {\cal M} }}

\newcommand{\nb}{\nabla}

\newcommand{\la}{\left\langle}
\newcommand{\ra}{\right\rangle}

\newcommand{\proof}{\emph{Proof. }}
\newcommand{\cvd}{\hfill$\square$ \bigskip}

\newcommand{\sss}{\mathbb S}
\newcommand{\qq}{\mathbb H}
\newcommand{\kk}{\mathbb K}
\newcommand{\pp}{\mathbb P}
\newcommand{\cc}{\mathbb C}
\newcommand{\cp}{\mathbb {CP}}

\newcommand{\sfe}{\mathbb S}

\newcommand{\tm}{T_{max}}

\newcommand{\aaa}{\left| A\right|^2}

\newcommand{\aaaq}{\left| A\right|^4}
\newcommand{\aao}{\accentset{\ \circ}{\left|A\right|}^2}

\newcommand{\hhh}{ H^2}

\newcommand{\nbh}{\left|\nabla H\right|^2}
\newcommand{\nba}{\left|\nabla A\right|^2}

\newcommand{\fse}{f_{\sigma,\eta}}
\newcommand{\fsep}{f_+}
\newcommand{\nbfse}{\left|\nb f_{\sigma,\eta}\right|^2}

\newcommand{\rr}{\bar r}
\newcommand{\dt}{\displaystyle{\frac{\partial}{\partial t}}}

\begin{document}

\def\qed{\hbox{\hskip 6pt\vrule width6pt height7pt
depth1pt  \hskip1pt}\bigskip}

%\setlength\parindent{0pt}
%\vspace*{1 cm}

\title{Cylindrical estimates for mean curvature flow of hypersurfaces in CROSSes}

\author{\sc G. Pipoli \footnote{
Institut Fourier, Universit\'e Joseph Fourier (Grenoble I), UMR 5582, CNRS-UJF, 38402, Saint-Martin-d'H\`eres, France. E-mail: giuseppe.pipoli@ujf-grenoble.fr}
 and C. Sinestrari
 \footnote{ (corresponding author) Dipartimento di Ingegneria Civile e Ingegneria Informatica, Universit\`a di Roma ``Tor Vergata'', Via Politecnico 1, 00133, Roma, Italy. E-mail: sinestra@mat.uniroma2.it}
 }
\date{}

\maketitle

\begin{abstract}
We consider the mean curvature flow of a closed hypersurface in the complex or quaternionic projective space. Under a suitable pinching assumption on the initial data, we prove apriori estimates on the principal curvatures which imply  
that the asymptotic profile near a singularity is either strictly convex or cylindrical. This result generalizes to a large class of symmetric ambient spaces the estimates obtained in previous works on the mean curvature flow of hypersurfaces in Euclidean space and in the sphere.

\end{abstract}

\medskip

\noindent {\bf MSC 2010 subject classification} 53C44, 35B40 \bigskip

\section{Introduction}\hspace{5 mm}  \setcounter{equation}{0}

We consider a family of closed hypersurfaces $\mm_t$ evolving by mean curvature flow in a Riemannian manifold. We are interested in the singular behaviour of the flow. More precisely, we want to describe a class of initial data which develops singularities whose blowup after rescaling is either strictly convex or cylindrical. 

When the ambient space is Euclidean, Huisken and Sinestrari \cite{HS2} showed that the above property holds provided the evolving surface is $2$-convex, i.e., the sum of the two smallest principal curvatures is strictly positive. This result is obtained by means of apriori estimates on suitable functions of the principal curvatures, called {\em cylindrical estimates}.  An alternative proof of this result has been later given by Haslhofer and Kleiner \cite{HK}. Recently, H.T. Nguyen \cite{Ng} has considered mean curvature flow in the sphere and proved cylindrical estimates under a pinching assumption which is related to $2$-convexity, but which is less restrictive in the regions with small curvature. In addition, cylindrical estimates have been obtained for suitable classes of nonlinear curvature flows by Andrews and Langford \cite{AL} in Euclidean space and, more recently, by Brendle and Huisken \cite{BH} in arbitrary Riemannian ambient spaces.

In this paper, we prove cylindrical estimates for the mean curvature flow when the ambient space is the projective space over either the complex field $\cc$ or the algebra of quaternions $\qq$, under a pinching assumption as in \cite{Ng}. More precisely, let $\kk\pp^n(4c)$ be the $n$-dimensional projective space with $\kk=\cc$ (resp. $\kk=\qq$) having real dimension $2n$ (resp. $4n$) and sectional curvature between $c$ and $4c$ for any given $c>0$. We set $m=2n-1$ (resp. $m=4n-1$) to denote the dimension of a hypersurface $\mm$ in $\kk\pp^n$. We denote by $\ld_1 \leq \dots \leq \ld_m$ the principal curvatures of $\mm$, by $H$ the mean curvature and by $|A|^2$ the squared norm of the second fundamental form. Our main result is the following.

\begin{Theorem}\label{main}
Let $n\geq 4$, and let $\mm_0$ be a closed real hypersurface of $\kk\pp^n(4c)$ that satisfies
\begin{equation}\label{pinching_hyp_2}
H>0, \qquad \aaa<\frac{1}{m-2}\hhh+4c.
\end{equation}
Then the mean curvature flow with initial data $\mm_0$ develops a singularity in finite time $\tm$ and condition \eqref{pinching_hyp_2} holds on $\mm_t$ for any $0\leq t<\tm$. Moreover, for every $\eta>0$ there exists a constant $C_{\eta}$ that depends only on $\eta$ and $\mm_0$ such that
\begin{equation}\label{cyl}
\left|\lambda_1\right|\leq\eta H \qquad\Rightarrow\qquad\left(\lambda_i-\lambda_j\right)^2\leq \Lambda\eta H^2+C_{\eta},\quad\forall i,j\geq 2,
\end{equation}
everywhere on $\mm_t$, for any $0\leq t<\tm$,
for a constant $\Lambda$ that depends only on the dimension.
\end{Theorem}

Taking into account the convexity estimate of \cite{HS1}, the above theorem shows that at those points where the curvature is large the local shape of the hypersurface is either uniformly convex or cylindrical. Our result, together with the one of \cite{Ng}, gives a fairly complete picture of the cylindrical estimates for the mean curvature flow in an ambient manifold which is a compact rank one symmetric space. 

It is interesting to compare the pinching assumption \eqref{pinching_hyp_2} with the hypothesis of $2$-convexity used in \cite{HS2,HK,AL,BH}, which means that $\ld_1+\ld_2>0$. It is easy to see that, at least at a pointwise level, the two properties are not comparable. On one hand, in fact, \eqref{pinching_hyp_2} is satisfied at any point where the curvatures are small enough, regardless of the sign of  $\ld_1+\ld_2$. On the other hand, \eqref{pinching_hyp_2} may be violated even at points where all $\ld_i$'s are positive but not enough close to each other, as simple explicit examples show. However, there is a relation between the two properties at the points with large curvature: as we show in \eqref{2conv}, at such points the pinching condition \eqref{pinching_hyp_2} implies $2$-convexity.

As in \cite{HS2,Ng}, using the cylindrical estimate it is possible to prove a pointwise bound on the gradient of the second fundamental form as follows.

\begin{Theorem}\label{gradient}
Under the hypotheses of the previous theorem, there is $C>0$ depending only on $\mm_0$ such that
$$
\nba\leq C(\aaaq+1)
$$
everywhere on $\mm_t$, for any $0\leq t<\tm$.
\end{Theorem}

The interest of the above results lies in possible applications to the construction of a generalized curvature flow after singularities by means of a surgery procedure. In fact, the estimates of the present paper and the ones of \cite{Ng}, together with the arguments of Section 8 in \cite{BH}, suggest that a flow with surgeries similar to the one of \cite{HS2} should exist also for the hypersurfaces in CROSSes satisfying the dimension restrictions and the pinching assumptions stated here. This would imply that any such hypersurface is diffeomorphic to a sphere $\sfe^m$ or to a finite connected sum of $\sfe^{m-1} \times \sfe^1$.

It should be remarked that the present analysis relies on the special structure of the ambient space considered. For a general nonsymmetric ambient manifold, the nonlinear flow recently introduced in \cite{BH} appears to be a more powerful tool for the study of $2$-convex hypersurfaces. Nevertheless, the pinching condition \eqref{pinching_hyp_2} retains an independent interest beside $2$-convexity, and it can lead to sharper results in symmetric spaces, as suggested by specific examples, see \S 3.2 in \cite{Ng}.

The methods of this paper are inspired by the above quoted papers \cite{HS2,Ng} and by Huisken's earlier work \cite{H3} on the mean curvature flow on the sphere. However, since the ambient spaces considered here are no longer space forms, the study presents additional difficulties, which lead to the dimension restriction $n \geq 4$ in our main theorem. We have also used a new argument for the lower bound on the polynomial $Z$ defined in \eqref{polZ}, see Lemma \ref{bound_Z}, which simplifies the corresponding step in \cite{Ng}, see Lemma 6.2 in that paper. In general, in our exposition we shall focus on the statements which require new arguments, and refer to the previous literature for the rest.

\section{Preliminaries}  \setcounter{equation}{0}

Let $F_0:\mm \to \left(\bmm,\bar g\right)$ be a smooth immersion of a closed connected real $m$-dimensional hypersurface into a Riemannian manifold and let $\nu$ be a normal unit vector field. In addition, let $g$ be the metric on $\mm$ induced by the immersion and let $A$ be the second fundamental form. We denote by $g_{ij}$ and $h_{ij}$ respectively the components of $g$ and $A$ in a given local coordinate system, and we call $\ld_1 \leq \dots \leq \ld_m$ the principal curvatures of $\mm$. Then the mean curvature $H$  associated with the immersion is $H=\ld_1+\dots+\ld_m$, and the squared norm of $A$ is given by $|A|^2=\ld_1^2+\dots+\ld^2_m$. In addition, we denote by $\aao=|A|^2-\frac 1m H^2$ the squared norm of the trace-free second fundamental form; then we have $\aao=\frac 1m \sum_{i<j}(\ld_i-\ld_j)^2 \geq 0$.

The evolution of $\mm_0 = F_0(\mm)$ by mean curvature flow is the one--parameter family of  immersions $F:\mm \times [0,\tm[ \,\to\bmm$ satisfying
 \begin{equation}\label{1.1}
\left\{\begin{array}{l}
\dt F(p,t)= -H\nu,
\qquad p \in \mm, \, t \geq 0, \medskip \\
F(\cdot,0)=F_0.
\end{array}\right.
\end{equation}
We denote by $\mm_t=F(\mm,t)$ the evolution of $\mm_0$ at time $t$.
It is well known that this problem has a unique smooth solution up to some maximal time $\tm \leq \infty$. Moreover, if $\tm$ is finite the curvature of $\mm_t$ necessarily becomes unbounded as $t \to \tm$ and we say that the flow develops a singularity. 

The ambient manifolds that we consider in this paper are $\bmm=\cp^n$, the complex projective space, and $\bmm=\qq\pp^n$, the quaternionic projective space. They belong to the class called Compact Rank One Symmetric Spaces (CROSSes).
An excellent introduction to these spaces can be found in chapter 3 of \cite{Be1}, here we recall the properties needed for our purposes. 

Let $\kk$ be either the field $\cc$ or the associative algebra $\qq$ and let $a$ be the real dimension of $\kk$, that is
$
a=2$ if $\kk=\cc$, while $a=4$ if $\kk=\qq$.
%$$
%a=\left\{\begin{array}{lr}
%2, & \text{if }\kk=\cc;\\
%4, & \text{if }\kk=\qq.
%\end{array}\right.
%$$   
We denote with $\mathbb{S}^n(c)$ the $n$-dimensional sphere with the canonical metric of constant curvature $c>0$. 
The action $T:\sss^{a-1}(1)\times\sss^{na+a-1}(c)\rightarrow\sss^{na+a-1}(c)$, $(\lambda,z)\mapsto\lambda z$ is by isometries and acts transitively on the fiber. $\kk\pp^n$ can be identified with $\sss^{na+a-1}/\sss^{a-1}$. The \emph{Hopf fibration} is the map $\pi:\sss^{na+a-1}(c)\rightarrow\kk\pp^n$, $z\mapsto\left[z\right]$, where $\left[z\right]$ is the class of $z$ under the action $T$. The Riemannian metric that we consider on $\kk\pp^n$ is the one induced from the metric of $\sss^{na+a-1}(c)$ such that $\pi$ becomes a Riemannian submersion. For $\kk=\cc$ this is the well-known Fubini-Study metric. 

\begin{Notation}
We denote by $\kk\pp^n(4c)$ the $\kk -$ projective space endowed with this metric. For simplicity of notation, we restrict ourselves to the case $c=1$ throughout the paper, and write $\kk\pp^n$ instead of $\kk\pp^n(4)$.  
\end{Notation}

The metric defined in this way has positive bounded sectional curvatures. In fact, if $X$ and $Y$ are two orthogonal unit vector tangent to $\kk\pp^n$, Theorem 3.30 of \cite{Be1} shows that $\overline {K}(X,Y)$, the sectional curvature of the tangent plane spanned by $X$ and $Y$, is
\begin{equation}\label{K_generale}
\overline {K}(X,Y)=1+3\left|pr_{Y\kk}X\right|^2,
\end{equation}
where $\left| \cdot \right|^2$ is the norm induced by the metric, and $pr_{Y\kk}X$ is the projection of $X$ on the tangent subspace $Y\kk$ of dimension $a$. 
It follows that $1\leq \overline {K}\leq 4$ and $\overline {K}=1$ (respectively $\overline {K}=4$) if and only if $X$ is orthogonal (respectively belongs) to $Y\kk$. 

We also recall that the CROSSes are Einstein manifolds. In the cases we are considering, the Einstein constant of $\kk\pp^n$ is
\begin{equation}\label{def_rr}
\rr=\left\{\begin{array}{rl}
2(n+1) =m+3 & \text{if }\kk=\cc, \medskip \\
4(n+2) = m+9 & \text{if }\kk=\qq.
\end{array}\right.
\end{equation}

The evolution equations of the main geometric quantities associated to a hypersurface evolving by mean curvature flow in a general Riemannian space have been derived in \cite{H2}. We recall here the equations satisfied by $H, \aaa$ and by the volume form $d\mu_t$. 

\begin{Lemma}\label{evoluzione_AH_hyp}
On a hypersurface evolving by mean curvature flow in a symmetric Einstein manifold we have
\begin{align}
&\nonumber 1)\qquad\dt H=\Delta  H+ H\left(\aaa+\rr \right),\\
%&\nonumber 1)\qquad \dt \tilde H=\Delta \tilde H+\tilde H\left(\aaa+\bar Ric(\nu,\nu)\right),\\
&\nonumber 2)\qquad \dt \hhh=\Delta \hhh-2\nbh+2\hhh\left(\aaa+\rr \right),\\
%\end{align}
%\begin{equation}
&\nonumber 3)\qquad\dt\aaa = \Delta\aaa-2\nba+2\aaa\left(\aaa+\rr \right) -4\left(h^{ij}h_j^{\phantom j p}\bar R_{pli}^{\phantom{pli}l}-h^{ij}h^{lp}\bar R_{pilj}\right),\\
%&\nonumber\phantom{3)\qquad\dt\aaa =} -4\left(h_{ij}h_j^{\phantom j p}\bar R_{pli}^{\phantom{pli}l}-h^{ij}h^{lp}\bar R_{pilj}\right),\\
%&\nonumber\phantom{3)\qquad\dt\aaa =} -4\left(h_{ij}h_j^{\phantom j p}\bar R_{pli}^{\phantom{pli}l}-h^{ij}h^{lp}\bar R_{pilj}\right),\\
&\nonumber 4)\qquad 
\dt d\mu_t=-\hhh d\mu_t,
\end{align}
where $\rr$ is the Einstein constant of the ambient manifold.
\end{Lemma}

\begin{Corollary}
If $\mm_0 \subset \kk\pp^n$ is a closed hypersurface with positive mean curvature, then its evolution by mean curvature flow becomes singular at a time $\tm <+\infty$. In addition, $\mm_t$ has positive mean curvature for all $t \in [0,\tm[\,$.
\end{Corollary}
\proof Let us set $H_{\min}(t)=\min_{\mm_t}H$. Using the evolution equation for $H$ in the previous lemma and the inequality $H^2 \leq m |A|^2$, we find
$$
\frac{d}{dt} H_{\min}(t) \geq \frac 1m H^3_{\min}(t).
$$ 
Since $H_{\min}(0)>0$ by our assumption, a standard comparison argument gives that $H_{\min}(t)$ remains positive and blows up in finite time. \cvd 

An important property in the analysis of singularities of the mean curvature flow of mean convex hypersurfaces are the convexity estimates proved in \cite{HS1}. They can be stated as follows.

\begin{Theorem} Let $F_0:\mm\to\bmm$ a smooth closed hypersurface immersion with non-negative mean curvature. Then, for any $\eta>0$ there is a $M_{\eta}$ such that
\begin{equation}\label{bound_l1}
\lambda_1\geq -\eta H-M_{\eta},
\end{equation}
everywhere on $\mm\times\left[0,\tm\right[$.
\end{Theorem}

The proof in \cite{HS1} is done in the Euclidean case but, as observed there, it is easily extended to the case of a general Riemannian ambient manifold. Finally, we recall the following inequality relating $\nba$ and $\nbh$, which is independent on the flow and is proved in Lemma 2.2 of \cite{H2}.

\begin{Lemma}\label{lemmaacod}
Let $\bmm$ be an Einstein manifold and $\mm$ a
 hypersurface of $\bmm$ of dimension $m$. Then at every point of $\mm$ we have
\begin{equation}
\nba \geq \frac{3}{m+2}\nbh.
\end{equation}
\end{Lemma}

\section{Preservation of pinching}\setcounter{equation}{0}\setcounter{Theorem}{0}

We begin the proof of Theorem \ref{main} by showing that our pinching inequality is preserved by the flow. We first observe that, since we assume the strict inequality \eqref{pinching_hyp_2}, by compactness we also have
\begin{equation}\label{base02_ep}
\aaa\leq a_{\ep}\hhh+b_{\ep}
\end{equation}
everywhere on $\mm_0$, where
$$
a_{\ep}=\frac{1}{m-2+\ep}, \qquad b_{\ep}=4(1-\ep),
$$
for some $\ep>0$ small enough depending only on $\mm_0$.

\begin{Proposition}\label{pinching_preserved_hyp}
Let $\mm_0$ be a closed hypersurface of $\kk\pp^n$, with $n\geq 4$. Then the pinching condition \eqref{base02_ep} is preserved by the mean curvature flow for any $\ep>0$ small enough.
\end{Proposition}
\proof  
Let us set $Q=\aaa-a_\ep \hhh-b_\ep$. Lemma \ref{evoluzione_AH_hyp} gives
\begin{equation}\label{evoluzione_Q_hyp}
\begin{array}{rcl}
\dt Q & =&  \Delta Q-2\left(\nba-a_\ep\nbh\right)+2\left(\aaa-a_\ep\hhh\right)\left(\aaa+\rr\right)\\
 & & -4\left(h^{ij}h_j^{\phantom j p}\bar R_{pli}^{\phantom{pli}l}-h^{ij}h^{lp}\bar R_{pilj}\right).
\end{array}
\end{equation}

By Lemma \ref{lemmaacod}, the contribution of the gradient terms in the above equation is non-positive. In addition, as in the proof of Proposition 3.6 in \cite{PS}, we can estimate
\begin{equation}\label{stimaI}
h^{ij}h_j^{\phantom j p}\bar R_{pli}^{\phantom{pli}l}-h^{ij}h^{lp}\bar R_{pilj} \geq m\aao.
\end{equation}
Thus we obtain
\begin{equation}\label{mp}
\dt Q  \leq  \Delta Q+2Q\left(\aaa+\rr\right)+2b_\ep\left(\aaa+\rr\right)  -4m\aao.
\end{equation}
Using the definition of $a_\ep,b_\ep$ and formula \eqref{def_rr}, we find
$$
2b_\ep -4m < - \frac 4{a_\ep}, \qquad 2\rr < \frac{4}{a_\ep},
$$
where we have also used the hypothesis $n \geq 4$, which means $m \geq 7$ and $m \geq 15$ in the cases $\kk=\cc$ and $\kk=\qq$ respectively. 
Therefore
$$
2b_\ep\left(\aaa+\rr\right)  -4m\aao =  ( 2b_\ep -4m)\aaa + 4 H^2 +2\rr b_\ep \leq  -\frac 4{a_\ep} Q. 
$$
Then we can apply the maximum principle to inequality \eqref{mp} and conclude that the inequality $Q \geq 0$ is preserved under the flow.
\cvd

Let us now define $W=\alpha\hhh+\beta$, where $\alpha,\beta$ are chosen such that 
\begin{equation}\label{alphabeta02}
\displaystyle{\frac{2}{m(m-2)}-\eta} <\alpha< \displaystyle{\frac{3}{m+2}-\frac{1}{m-1}-\eta,} \end{equation}
\begin{equation}\label{beta}
\beta=b_\ep=4(1-\ep).
\end{equation}
In addition, we introduce the auxiliary function $$\displaystyle{\fse:=\frac{\aaa-(\frac{1}{m-1}+\eta)\hhh}{W^{1-\sigma}}},$$ where $\sigma$ and $\eta$ are two positive constants small enough. For simplicity we write $f_0=f_{0,\eta}$. From our choice of $\alpha,\beta$ and Proposition \ref{pinching_preserved_hyp} it is easily checked that
\begin{equation}
\label{fsigma}
f_{\sigma,\eta} \leq W^\sigma.
\end{equation}
We now derive a differential inequality satisfied by $\fse$.
In the rest of the paper, we denote by $C_1,C_2,\dots$ positive constants depending on the dimension and the initial data, but not on $\eta,\sigma$ provided they are small enough.

\begin{Proposition}\label{evoluzione_fse}
There exists a constant $\sigma_1 \in \,]0,1[\,$, depending only on $\mm_0$ such that at the points where $\fse>0$  
\begin{equation}\label{dt_fse}
\begin{array}{rcl}
\dt \fse&\leq&\displaystyle{ \Delta \fse+\frac{2\alpha (1-\sigma)}{W}\la\nb \fse,\nb \left|H\right|^2\ra-2C_1W^{\sigma-1}\left|\nb H\right|^2}\\
 & & \displaystyle{+2\sigma(\aaa+\rr) \fse-2C_2\fse,}
\end{array}
\end{equation}
holds for  all $0\leq\sigma\leq\sigma_1$, for suitable constants $C_1,C_2>0$.
\end{Proposition}
\proof Making similar calculations to the proof of Proposition 4.1 of \cite{PS} we have
\begin{eqnarray*}
\dt\fse & \leq &  \displaystyle{\Delta\fse+\frac{2\alpha(1-\sigma)}{W}\la\nb\fse,\nb\hhh\ra -2W^{\sigma-1}\nba}\\
 & & \displaystyle{+2W^{\sigma-1}\left(\frac{1}{m-1}+\eta+f_0(1-\sigma)\alpha \right) \nbh}\\
 & & +2\beta(1-\sigma)\frac{\fse}{W}\left(\aaa+\rr\right)+2\sigma\fse(\aaa+\rr)-4mW^{\sigma-1}\aao.
\end{eqnarray*}
With the choice \eqref{alphabeta02} of $\alpha$ and $\beta$ we have $f_0<1$. Hence by Lemma \ref{lemmaacod}
\begin{eqnarray*}
\lefteqn{-\left|\nb A\right|^2 +\displaystyle{\left (\frac{1}{m-1}+\eta+f_0(1-\sigma)\alpha\right)\left|\nb H\right|^2}}\\
 &\leq& \displaystyle{\left(\frac{1}{m-1}+\eta+\alpha\right)\left|\nb H\right|^2-\left|\nb A\right|^2}\\
 & \leq &\displaystyle{\left(\alpha+\frac{1}{m-1}+\eta-\frac{3}{m+2}\right)\left|\nb H\right|^2 }\\
 & = & \displaystyle{-C_1\left|\nb H\right|^2,}
\end{eqnarray*}
with $C_1>0$ by  \eqref{alphabeta02}. We now consider the reaction terms. Using conditions \eqref{base02_ep}, \eqref{alphabeta02}, \eqref{beta} and the fact that $\fse\leq W^{\sigma}$, we find, for both choices of $\kk$,
\begin{eqnarray*}
\lefteqn{2\beta(1-\sigma)\frac{\fse}{W}\left(\aaa+\rr\right)-4mW^{\sigma-1}\aao} \\
 & \leq & 2\beta(1-\sigma)\frac{\fse}{W}\left(\aaa+\rr\right)-4m(\fse+\eta\hhh W^{\sigma-1}) \\
  & \leq & 2\beta\frac{\fse}{W}\left(\aaa+\rr\right)-4m\frac{\fse}{W}(W+\eta\hhh) \\
 &\leq & 2\frac{\fse}{W}\left[\beta\left(\frac{\hhh}{m-2+\ep}+4(1-\ep)+\rr\right)\right]\\
 & & -4m\frac{\fse}{W}\left[\left(\frac{2}{m(m-2)}-\eta\right)\hhh+\beta +\eta \hhh \right]\\
 &\leq & 2\frac{\fse}{W}\left[ \left( \frac{\beta}{m-2} -\frac{4}{m-2} \right)\hhh +\left( 4(1-\ep)+\rr -2m\right) \beta \right] \\
 &\leq & -2\frac{\fse}{W}\left[ \frac{4 \ep}{m-2} \hhh +4\ep \beta \right] \\
  & \leq  & -4\ep \fse.
\end{eqnarray*}
Collecting the above inequalities, the assertion follows. \cvd

\section{Cylindrical estimates}\setcounter{equation}{0}\setcounter{Theorem}{0}

In this section we prove the following result, from which the main assertions of our paper follow easily.

\begin{Theorem}\label{cyl_th}
Let $n\geq 4$ and  $\mm_0$ a closed hypersurface of $\kk\pp^n$. If $\mm_0$ satisfies conditions \eqref{pinching_hyp_2}, then for any $\eta>0$ there exists a constant $C'_{\eta}$ depending on $\eta$ and the initial data such that
$$
\aaa-\frac{1}{m-1}\hhh\leq\eta\hhh+C'_{\eta},
$$
for every time $t\in\left[0,\tm\right[$.
\end{Theorem}

To prove the above statement, we will show that the function $\fse$ introduced in the previous section is bounded. To this purpose, we follow the procedure first introduced in \cite{H1}, which relies on integral estimates and Stampacchia iteration. We recall Simons' identity, see e.g. \cite[Lemma 2.1]{H2},
\begin{eqnarray*}
\Delta\aaa  & = &  2\la h_{ij},\nb_i\nb_j H\ra+2\nba+2Z+2\left(Hh^{ij}\bar R_{0i0j}-\aaa\bar R_{0l0}^{\phantom{0l0}l}\right),
\\
&& +4\left(h^{ij}h_j^{\phantom j p}\bar R_{pli}^{\phantom{pli}l}-h^{ij}h^{lp}\bar R_{pilj}\right),
\end{eqnarray*}
where 
\begin{equation}\label{polZ}
\displaystyle Z=\left(\sum_{i=1}^m \ld_i \right) \left(\sum_{i=1}^m \ld_i^3 \right)-\aaaq.
\end{equation}
By the boundedness of the curvature of the ambient space and the pinching condition \eqref{base02_ep}, we have that
$$
Hh^{ij}\bar R_{0i0j}-\aaa\bar R_{0l0}^{\phantom{0l0}l} \geq-C_3W,
$$
where $C_3$ is some positive constant. Using this and \eqref{stimaI}, we find
\begin{equation}\label{Delta_A}
\displaystyle{\Delta \aaa \geq 2\la h_{ij},\nb_i\nb_j H\ra+2\nba+2Z-2C_3 W.}
\end{equation}
Unlike the case considered in \cite{H1}, the polynomial $Z$ can change sign under our assumptions. However, combining the pinching condition with the convexity estimate, we can prove the following inequality which shows that the negative part of $Z$ is of lower order.

\begin{Lemma}\label{bound_Z}
Assuming the pinching condition \eqref{base02_ep} and $H>0$, there is a constant $\gamma$, depending only on $n$ and $\mm_0$, such that for any $\eta>0$ there exists a $K_{\eta}$ such that
\begin{equation}\label{z1}
Z\geq \gamma\hhh\left(\aaa-\frac{1}{m-1}\hhh-\eta\hhh\right)-K_{\eta}\left(H^3+1\right).
\end{equation}
Therefore, there are also constants $\gamma',K'_{\eta}$ such that
\begin{equation}\label{z2}
Z\geq \gamma' W^2 f_0-K'_{\eta}\left(H^3+1\right).
\end{equation}
\end{Lemma}
\proof
We first show that inequality \eqref{base02_ep} implies $2$-convexity at the points where $H$ is large. We have the identity
\begin{eqnarray}
\aaa - \frac{1}{m-2}H^2 & = & \ld_1^2+\ld_2^2 + \frac{1}{m-2} (\ld_1+\ld_2)^2 \nonumber \\
& &- \frac{2}{m-2}(\ld_1+\ld_2)H + \frac{1}{m-2} \sum_{3 \leq i <j \leq m}(\ld_i-\ld_j)^2.
\end{eqnarray}
Neglecting the positive terms in the above equality, and recalling Proposition \ref{pinching_preserved_hyp}, we obtain
$$
(\ld_1+\ld_2)H \geq -\frac{m-2}{2}\left(\aaa - \frac{1}{m-2}H^2 \right) \geq \frac{\ep}{2(m-2+\ep)}H^2 -2(m-2).
$$
It follows that at any point where $\displaystyle H^2 \geq \frac8\ep{(m-2)(m-2+\ep)}$ we have
\begin{equation}\label{2conv}
\ld_1+\ld_2 \geq \frac{\ep}{2(m-2+\ep)} \left(H -\frac{H}{2} \right) =  \frac{\ep}{4(m-2+\ep)} H.
\end{equation}
Once we have this estimate, we can argue exactly as in Lemma 5.2 in \cite{HS2} to conclude that at such points the assertion holds.

It remains to consider the part of the hypersurface where $H^2 < 8(m-2)(m-2+\ep)/\ep$. This case is trivial, because on this set $|A|^2$, and therefore $|Z|$, is uniformly bounded by some constant only depending on $m$ and the initial data, thanks to \eqref{base02_ep}. Finally, estimate \eqref{z2} 
follows easily from \eqref{z1} and the definition of $W$. \cvd

{\em Proof of Theorem \ref{cyl_th}.} Once the above estimates have been established, the result can be obtained from the method of proof as in \cite{HS2,Ng} with easy modifications. We sketch the main steps for the reader's convenience. We first derive an inequality for $\Delta \fse$ by using estimate \eqref{Delta_A}
\begin{equation}\label{Delta_fse02}
\begin{array}{rcl}
\Delta \fse & \geq &\displaystyle{ 2 W^{\sigma-1}\la h_{ij}^{\eta},\nb_i\nb_j H\ra+2ZW^{\sigma-1}-2C_3W^{\sigma}}\\
 & &\displaystyle{ -2(1-\sigma)\frac{H\fse}{W}\Delta H-4\alpha(1-\sigma)\frac{H}{W}\la\nb_i H,\nb_i\fse\ra,}
\end{array}\end{equation}
where $\displaystyle{h_{ij}^{\eta}=h_{ij}-\left(\frac{1}{m-1}+\eta\right)H g_{ij}}$. We now denote by $f_+$ the positive part of $\fse$ and proceed to  estimate suitable $L^p$ norms of $f_+$. To this aim, we multiply inequality \eqref{Delta_fse02} by $\fsep^{p-1}$, for $p>1$ large enough, 
and integrate over $\mm_t$. A long computation, similar to  \cite[Lemma 5.5]{HS2} or \cite[Proposition 6.7]{Ng}, allows to find positive constants $C_4$ and $C_5$ such that for every $\delta>0$ the following holds:
\begin{eqnarray} \nonumber
\frac{1}{C_4}\int_{\mm_t} W f_+^p d\mu & \leq & (1+\delta )\int_{\mm_t} f_+^{p-1}W^{\sigma-1}\nbh d\mu \\
& &  +\frac{1}{\delta}\int_{\mm_t} f_+^{p-2}\nbfse d\mu+C_5. \label{poincare}
\end{eqnarray}
A crucial step in deriving the above estimate is provided by Lemma \ref{bound_Z}, which allows to relate the left-hand side of \eqref{poincare} to the term containing the factor $Z$ which comes from the integration of \eqref{Delta_fse02}.

Then we can estimate the $L^p$ norm of $\fsep$. We multiply the inequality in Theorem \ref{evoluzione_fse} by $p f_+^{p-1}$ and integrate over $\mm_t$. Inequality \eqref{poincare} allows to absorb the positive zero order term in the equation with the negative gradient terms, for suitable choices of $\sigma, p$ and $\delta$. More precisely, with computations similar to \cite[Lemma 5.6]{HS2} or \cite[Proposition 6.8]{Ng}, we obtain that there are positive constants $C_6$, $C_7$ and $C_8$ depending only on $\mm_0$ such that, for all $p\geq C_7$ and $\sigma\leq \frac{C_8}{\sqrt{p}}$, we have
$$
\left(\int_{\mm_t} f_+^pd\mu\right)^{\frac{1}{p}}\leq C_6, \qquad  \forall t \in [0,T_{max}[\,.
$$

Once this $L^p$ estimate is established, a Stampacchia iteration procedure as in \cite[Theorem 5.1]{H1} shows that $\fse$ is uniformly bounded from above for an appropriate choice of $\sigma$. By the definition of $\fse$, this property implies Theorem \ref{cyl_th}. \cvd

The main results stated in the introduction are an easy consequence of Theorem \ref{cyl_th}. \bigskip

\noindent{\em Proof of Theorem \ref{main}}.  The following identity holds in general
\begin{equation}\label{2001}
\aaa-\frac{1}{m-1}\hhh = \frac{1}{m-1}\left(\sum_{1<i<j}\left(\lambda_i-\lambda_j\right)^2+\lambda_1(m\lambda_1-2H)\right),
\end{equation}
where $\lambda_1\leq\dots\leq\lambda_m$ are the principal curvatures of $\mm_t$. Then for every $i>j>1$ we have
$$
\left(\lambda_i-\lambda_j\right)^2\leq(m-1)\left(\aaa-\frac{1}{m-1}\hhh\right)-\lambda_1(m\lambda_1-2H).
$$
Using Theorem \ref{cyl_th} and $\left|\lambda_1\right|<\eta H$ we have
$$
\left(\lambda_i-\lambda_j\right)^2\leq(m-1)\left(\eta\hhh+C'_{\eta}\right)+\eta(m\eta+2)\hhh,
$$
which gives the assertion for suitable constants $\Lambda$ and $C_{\eta}$. \cvd

\noindent{\em Proof of Theorem \ref{gradient}}.  The estimate of Theorem \ref{cyl_th} allows us to consider a function of the form
$$
g:=\frac{|\nabla A|^2}{ ( c_1 H^2 -|A|^2+K_1)(c_2 H^2-|A|^2+K_2)}
$$
for suitable constants
$
c_i > \frac{1}{m-1}$ and $K_i>0$. As in Theorem 6.1 of \cite{HS2} or Theorem 7.1 of \cite{Ng}, we can compute the evolution equation for $g$ and apply the maximum principle to find an upper bound on this function, which gives the desired estimate on $|\nabla A|^2$. \cvd

\noindent {\bf Acknowledgments} 
The results of this paper are part of Giuseppe Pipoli's PhD thesis, written at the Department of Mathematics, University ``Sapienza'' of Rome. 
Giuseppe Pipoli was partially supported by PRIN07 ``Geometria Riemanniana e strutture differenziabili'' of MIUR (Italy) and Progetto universitario Univ. La Sapienza ``Geometria differenziale -- Applicazioni''. Carlo Sinestrari was partially supported by FIRB--IDEAS project ``Analysis and beyond'' and by the group GNAMPA of INdAM (Istituto Nazionale di Alta Matematica).

\bigskip

\end{document}